\documentstyle[leqno,11pt]{article}

\newtheorem{Thm}{Theorem}[section]
\newtheorem{Lemma}[Thm]{Lemma}
\newtheorem{Prop}[Thm]{Proposition}
\newtheorem{Cor}[Thm]{Corollary}

\newcommand{\surj}{\,\longrightarrow\hspace{-12pt}
                     \longrightarrow\,\,}

\newcommand{\pf}{\ \\ \noindent {\bf Proof.\ \ }}
\newcommand{\qed}{\ $\displaystyle\Box$\\ \ \par}

\newcommand{\rk}{\refstepcounter{Thm}
\ \\ \noindent {\bf Remark \theThm }\ \ }

\renewcommand{\H}{{\cal H}}
\newcommand{\C}{{\bf C}}
\newcommand{\Z}{{\bf Z}}
\newcommand{\X}{{\cal X}}
\newcommand{\U}{{\cal U}}
\newcommand{\OO}{{\cal O}}
\newcommand{\res}{\mbox{\rm res}}
\newcommand{\Res}{\mbox{\rm Res}}

\title{\bf Riemann reciprocity in higher dimensions}

\author{{\sc Yakov Karpishpan}\\ \\ \small
Department of Mathematics,  Yeshiva University\\
\small 500 W 185 St., New York, NY 10033, USA}

\date{\em Elul 5754}

\begin{document}

\maketitle

\begin{abstract}
The reciprocity law for abelian
differentials of first and second kind is generalized to
higher-dimensional
varieties. It is shown that $H^1(V)$ of a polarized variety $V$ is
encoded in the Laurent data  along a curve germ in $V$, with the
polarization form on $H^1(V)$ corresponding to the {\em
one-dimensional} residue pairing. This associates an {\em extended
abelian variety} to $V$; if $V$ is an abelian variety itself, our
construction ``extends" it, even when $V$ is not a Jacobian.
\end{abstract}

\section{Introduction}
The reciprocity law, or bilinear relation, for abelian
differentials of the first and second kind is classically
formulated as follows. Let $X$ be a compact Riemann surface of
genus $g$.
Suppose $\omega$ is a holomorphic one-form (an abelian differential
of the first kind) and $\eta$ is a meromorphic one-form with
exactly one pole, necessarily of order $\geq 2$, at some point
$p\in X$ (an abelian differential of the second kind). Let
$A_1,\ldots,A_g, B_1\ldots, B_g$ (resp. $A'_1,\ldots,A'_g,
B'_1,\ldots, B'_g$) denote the periods of
$\omega$ (resp.  $\eta$) over a standard symplectic basis for
$H_1(X,\Z)$. Then
\begin{equation}
\label{classic}
\sum_{j=1}^{g}A_jB'_j - A'_jB_j = 2\pi i \,\res_p\ (f\eta)\ ,
\end{equation}
where $f$ is a holomorphic function defined near $p$, with
$df=\omega$.

The left-hand side of (\ref{classic}) can be more easily recognized
for what it is when written as a matrix product
$$
(A_1 \ldots B_g)\left(
\begin{array}{cc}
0 & I_g\\
-I_g & 0
\end{array}
\right)
\left(\begin{array}{c}
A'_1\\
\vdots\\
B'_g
\end{array}\right)\ .
$$
The forms $\omega$ and $\eta$ represent cohomology classes in
$H^1(X,\C)$ determined by their vectors of periods. The
polarization form $Q$ is defined as the  dual of the
intersection form on $H_1(X,\Z)$, and so the matrix product above
equals $Q([\omega],[\eta])$.

As to the right-hand side of (\ref{classic}), let $g$ be a
meromorphic function defined near $p$ with $dg=\eta$. Writing
$<f,g>$ for the residue pairing $\res_p (fdg)$, one has the
following version of
(\ref{classic}):
\begin{equation}
\label{modern}
Q([\omega],[\eta])=2\pi i\,<f,g>\ .
\end{equation}

It is this statement that we wish to generalize in this paper for
varieties of dimension greater than one. The pointed Riemann
surface is replaced by a smooth complex projective variety $V$ with
an irreducible ample divisor $D$ and a regular point $p$ on $D$ (we
do not assume that $D$ is smooth).
We also choose a smooth curve germ $\X$ on $V$ through $p$,
transversal to $D$.
The Chern class of $D$ defines  a polarization form $Q$ on
$H^1(V,\C)$. The claim is that $Q$ is again expressed via the {\em
one-dimensional} residue pairing. Specifically, let $\omega$ be a
holomorphic one-form (a simple abelian differential of the first
kind) on $V$, and let $\eta$ be a {\em closed} meromorphic one-form
on $V$ with poles only along $D$ (a simple abelian differential of
the second kind).
Such forms again represent cohomology classes in $H^1(V,\C)$
(Proposition \ref{V-D:V}), and
in  Theorem \ref{main} we establish a similar relation
\begin{equation}
\label{our}
Q([\omega],[\eta])=(-1)^{n-1}2\pi i\,<f,g>\ ,
\end{equation}
where now $f$ and $g$ are, respectively, holomorphic and
meromorphic functions on $\X$, whose differentials are the
one-forms $\omega$ and  $\eta$ {\em pulled back to} $\X$:
$$
df=\omega|_{\X}\ ,\hspace{1in}dg=\eta|_{\X}\ ,
$$
and $<f,g>=\res_pfdg$ denotes the residue pairing on $\X$ at $p$.

Returning to the one-dimensional case, with additional notation and
terminology a more complete statement is possible.  Let
$\H=\C((z))$ (the field of formal Laurent power series),
$\H_+=\C[[z]]$, and $\H'=\C((z))/\C$. The pairing
$<f,g>=\res_{z=0}fdg$ on $\H$ is skew-symmetric; it is
non-degenerate on $\H'$.
Choosing a formal local parameter
$u:\hat{\OO}_{X,p}\stackrel{\cong}{\longrightarrow}\H_+$
on $X$ at $p$ defines maps
$$
\Gamma(X,\OO_X(*p))\longrightarrow \H'\ \ \ \mbox{and}\ \ \
\Gamma(X,\Omega^1_X(*p))\longrightarrow \H dz\ ,
$$
both of which will also be denoted $u$.

Put $K_0=u(\Gamma(X,\OO_X(*p)))$ and
$\Omega=\{f\in\H'\,|\,df\in u(\Gamma(X,\Omega^1_X(*p)))\}$. Then
$K_0\subset \Omega$ are each other's annihilators in $\H'$, i.e.
$\Omega=K_0^{\perp}$ and $<\ ,\ >$ induces a perfect pairing on
$K_0^{\perp}/K_0$. The reciprocity law now says that there is a
symplectic isomorphism
\begin{equation}
\label{basic:iso}
(H^1(X,\C),\ Q)\longrightarrow (K_0^{\perp}/K_0,\ 2\pi i\,<\ ,\ >)\
{}.
\end{equation}

Noting that $(H^1(X,\C),Q)$ is a polarized Hodge structure,
we may go further and  identify all of its components in terms of
the Laurent data.
Thus $H^{1,0}=F^1H^1(X,\C)$ is simply $K_0^{\perp}\cap\H'_+$.
As to the integral structure, put
$K=\{f\in\H'\,|\,e^f\in u(\Gamma(X-\{p\},\OO^*))\}$. Then
$K_0\subset K\subset K_0^{\perp}$ and $\Lambda:=K/K_0$ is
isomorphic to $H^1(X,\Z)$.
Finally, let $U\subset K_0^{\perp}/K_0$ be the image of $H^{0,1}$
under (\ref{basic:iso}) and, in turn, let $Z$ denote the preimage
of $U$ under the projection
$K_0^{\perp}\rightarrow K_0^{\perp}/K_0$.
Then $Z$ is a maximal isotropic subspace of $\H'$.

These results, most of which may already be found in \cite{SW},
led Arbarello and De Concini \cite{AD} to codify triples
$(Z,K_0,\Lambda)$ with properties as above  under the name of
{\em extended abelian varieties}.
They are also called {\em extended Hodge structures of
weight one} in \cite{K}.
Evidently, to each such there corresponds a unique polarized Hodge
structure of weight one, although going back there are infinitely
many choices (see (\ref{choices})).
Thus one has a ``de-extension" map with infinite fibers

\begin{equation}
\label{de-extend}
\left\{\mbox{extended abelian varieties}\right\}\longrightarrow
\left\{\mbox{abelian varieties}\right\}\ .
\end{equation}

In this paper it is shown that all of the above generalizes to
higher-dimen\-sion\-al varieties.
We again meet the subspaces  $K_0\subset\Omega$ of $\H'$ with
$\Omega=K_0^{\perp}$, and a reciprocity law analogous to
(\ref{basic:iso}) holds (Theorem \ref{main:two}):
there is a symplectic isomorphism
\begin{equation}
\label{basic:iso:two}
(H^1(V,\C),\ Q)\longrightarrow (K_0^{\perp}/K_0,\ (-1)^{n-1}2\pi
i<\ ,\
>)\ .
\end{equation}
We also construct the remaining components of an
extended abelian variety whose associated Hodge structure is that
of $H^1(V,\C)$ for any variety $V$ with a divisor $D$ as above
(Theorem \ref{EHS}).
In particular, this applies when $V$ is an abelian variety and $D$
is its theta divisor. Thus we obtain a method of ``inverting" the
``de-extension" map (\ref{de-extend}), or ``extending" a Hodge
structure of weight one, even if the latter did not come from
geometry.
This may be useful in approaching the Schottky problem, which was,
in fact, our original motivation.

We leave with a question suggested by the present work: does this
generalize for differential forms of higher degree, and is there a
good notion of an extended
Hodge structure of weight higher than one?

\noindent{\bf Acknowledgements}. During work on this paper the
author had helpful conversations with Enrico Arbarello and Roberto
Silvotti. This project has also been influenced by the attempts of
M. Sato and his co-workers at a higher-dimensional generalization
of Krichever's theory, sketched in \cite{Sa}.

\section{$H^1$ of a polarized variety}
\label{H1}
We will  be studying a smooth complex projective variety $V$ with
an irreducible ample divisor $D$, not necessarily smooth. By
Grothendieck's Algebraic De Rham Theorem \cite{Gro},
\begin{equation}
H^k(V-D,\C)\cong\frac{\Gamma(\tilde{\Omega}_V^k(*D))}
{d\Gamma(\Omega_V^{k-1}(*D))}\ ,
\end{equation}
where the tilde in $\tilde{\Omega}_V^k(*D)$ denotes  the subsheaf
of $d$-closed forms in $\Omega_V^k(*D)$.

\begin{Prop}
\label{V-D:V}With the above assumptions,
$$
H^1(V,\C)\cong H^1(V-D,\C)\cong
\frac{\Gamma(\tilde{\Omega}_V^1(*D))}
{d\Gamma(\OO_V(*D))}\ .
$$
\end{Prop}

\pf
First, assume $D$ is smooth. Then there is a sequence of sheaf
complexes
$0\rightarrow\Omega_V^{\bullet}\rightarrow
\Omega_V^{\bullet}(\log
D)\rightarrow\Omega_D^{\bullet-1}\rightarrow 0$ inducing
$$
0\longrightarrow H^1(V)\longrightarrow
H^1(V-D)\stackrel{P.R.}{\longrightarrow}
H^0(D)\stackrel{\gamma}{\longrightarrow}
H^2(V)\longrightarrow\ldots
$$
Here $P.R.$ stands for the Poincar\'e Residue and $\gamma$ denotes
the Gysin map, which is obtained by applying Poincar\'e duality on
the source and the target of the map
$H_{2n-2}(D)\rightarrow H_{2n-2}(V)$, where $n=\dim_{\C}V$. The
latter map is injective: the image of a suitable generator of
$H_{2n-2}(D)(\cong\C)$ is the fundamental class of $D$ in
$H_{2n-2}(V)$, which cannot be zero by the irreducibility of $D$.
Hence $\gamma$ is injective too, which proves the assertion.

In the general case, when $D$ is not smooth, the above arguement
does not work. Instead, consider the spectral sequence
$$
E_2^{p,q}=H^p(V,\H^q\Omega_V^{\bullet}(*D))\Rightarrow H^*(V-D,\C)\
{}.
$$
It yields the very same exact sequence
$$
0\longrightarrow H^1(V)\longrightarrow
H^1(V-D)\stackrel{R}{\longrightarrow}
H^0(D)\stackrel{\nu}{\longrightarrow} H^2(V)\longrightarrow\ldots
$$
where now $R$ is the cohomological residue map induced by the sheaf
isomorphism $\H^1\Omega_V^{\bullet}(*D)\cong\C_D$, and $\nu$ is the
map sending $1_D$ to (the Poincar\'e dual of) the fundamental class
of $D$ in $H^2(V)$ (see \cite{GH}, p. 458). Again, the map $\nu$
must be injective, i.e. $H^1(V)\cong H^1(V-D)$ in all cases. \qed

We note that the isomorphism of the above Proposition makes the
Hodge filtration  on $H^1(V,\C)$ obvious:
$$
F^1H^1(V,\C)=\Gamma(V,\Omega_V^1)\ .
$$
There is also a corresponding isomorphism of quotients
\begin{equation}
\frac{\Gamma(V,\tilde{\Omega}^1_V(*D))}
{\Gamma(V,\Omega^1_V) + d\Gamma(V,\OO_V(*D))}
\stackrel{\cong}{\longrightarrow}H^1(V,\OO_V)\ .
\end{equation}
For future use, we will need to define it explicitly. Let
$\U=\{U_{\alpha}\}$ be an affine open cover of $V$ or an acyclic refinement of
such; thus $\check{H}^1(\U,\OO_V)\cong H^1(V,\OO_V)$.
The \v{C}ech cochain
$\{g_{\alpha\beta}\in\Gamma(U_{\alpha}\cap U_{\beta},\OO_V)\}$
giving the image of $\eta\in\Gamma(V,\tilde{\Omega}^1_V(*D))$ in
$H^1(V,\OO_V)$ can be described as follows:
$$
\{g_{\alpha\beta}\}=\check{\delta}\{\mu_{\alpha}\}\ ,
$$
where $\mu_{\alpha}$ is a meromorphic function on $U_{\alpha}$ such
that
$\eta|_{U_{\alpha}}-d\mu_{\alpha}$ is a holomorphic one-form on
$U_{\alpha}$.
In particular, when $U_{\alpha}$ is simply connected, we may assume
that $\eta|_{U_{\alpha}}=d\mu_{\alpha}$, and if $\eta$ is already
holomorphic over $U_{\alpha}$, then
we may take $\mu_{\alpha}=0$.

\section{Polarization and residues}
\label{polarization}
The Hodge structure on $H^1(V,\C)$ is polarized by the divisor $D$.
Concretely, the cup product
\begin{equation}
\label{QQ}
Q:H^1(V,\C)\stackrel{\smile}{\otimes} H^1(V,\C)
\stackrel{\smile c_{[D]}^{n-1}}{\longrightarrow}
H^{2n}(V,\C)
\stackrel{\int}{\longrightarrow}\C
\end{equation}
gives an integrally-defined perfect pairing. Here $c_{[D]}\in
H^2(V,\Z)$ stands for the first Chern class of $\OO_V(D)$,
and $\int$ denotes the topological trace, i.e. the map
$H^{2n}_{DR}(V,\C)\stackrel{\sim}{\rightarrow}\C$ obtained by
integrating $C^{\infty}$ complex-valued $2n$-forms over $V$.

Thinking of $c_{[D]}$ as an element of
$H^1(V,\Omega_V^1)$, we will generally prefer a more algebraic
version
\begin{equation}
\label{Q}
Q:H^0(V,\Omega_V^1)\stackrel{\smile}{\otimes} H^1(V,\OO_V)
\stackrel{\smile c_{[D]}^{n-1}}{\longrightarrow}
H^n(V,\Omega_V^n)\cong H^{2n}(V,\C)
\stackrel{\int}{\longrightarrow}\C
\end{equation}
\rk
No notational distinction is made between the classes $[\omega]$,
$[\eta]$ in
$H^1(V,\C)$ and in $H^0(V,\Omega_V^1)$ or $H^1(V,\OO_V)$. Likewise,
both
pairings above are denoted $Q$. This should cause no confusion,
since the
pairings are compatible with the natural maps
$H^0(V,\Omega_V^1)\hookrightarrow
H^1(V,\C)$ and $H^1(V,\C)\surj H^1(V,\OO_V)$.

Even more
algebraically, $\int=(2\pi i)^n tr$, where $tr$ is
Grothendieck's trace isomorphism (cf. \cite{D}, (2.2)). As the
latter is the least explicit part in the definition of $Q$, let us
elaborate on it.

First recall from \cite{L} and \cite{GH} a construction of the
local trace
$$
H^n_{\{p\}}(V,\Omega_V^n)
\stackrel{tr_p}{\longrightarrow}\C\ .
$$
Let $U$ be a polycylindrical neighborhood of $p$ in $V$. Then
$$
H^n_{\{p\}}(V,\Omega_V^n)\cong H^n_{\{p\}}(U,\Omega_U^n)\ ,
$$
by excision. And $U$ being Stein (and $n\geq 1$) implies
$$
H^n_{\{p\}}(U,\Omega_U^n)\cong H^{n-1}(U-\{p\},\Omega_U^n)\ .
$$
These isomorphisms are composed with the residue morphism
$$
\Res:
H^{n-1}(U-\{p\},\Omega_U^n)\stackrel{\cong}{\longrightarrow}\C\ .
$$

To define $\Res$, let
$t_1,\ldots,t_n$ be a coordinate system on $U$ centered at $p$.
Thus $D_j=\{t_j=0\}$ are $n$ smooth hypersufaces in $U$
intersecting normally at $p$.
We assume that $D_1=D\cap U$, although this will not be necessary
until the next section. Put $U_j=U-D_j$, and let $D^+=D_1+\ldots
+D_n$. Evidently, the $U_j$'s form an acyclic open cover of
$U-\{p\}$.
Then a class in $H^{n-1}(U-\{p\},\Omega_U^n)$ is represented by a
holomorphic form on $U_1\cap\ldots \cap U_n=U-D^+$.
In fact, we may assume that this form is a restricton to
$U_1\cap\ldots \cap U_n$  of a meromorphic form
$\psi\in\Gamma(U,\Omega_U^n(*D^+))$. The form $\psi$ can be
expanded in a Laurent series:
$$
\psi=\sum a_{k_1,\ldots, k_n}t^{k_1}\cdots
t^{k_n}dt_1\wedge\ldots\wedge dt_n\ .
$$
Finally,
$$
\Res \,\psi\stackrel{\rm def}{=} a_{-1,\ldots,-1}\ .
$$
It is a basic fact of residue theory that $\Res \,\psi$ is
independent of the parameter system $(t_1,\ldots,t_n)$.

We note one property of
$$
\Res: \Gamma(U,\Omega_U^n(*D^+))\longrightarrow \C
$$
for future reference. Assume $\psi\in\Gamma(U,\Omega_U(*D^+))$ has
only {\em simple} poles along $D_2,\ldots,D_n$, i.e. can be written
as
$$
\psi=\frac{h\,dt_1\wedge\ldots\wedge dt_n}{t_1^m\,t_2\cdots
t_n}\ ,
$$
with $h$ a holomorphic function.
Then
\begin{equation}
\label{Res:res}
\Res\,\psi = \res_0\, \frac{\tilde{h}\,dt_1}{t_1^m}\ ,
\end{equation}
where $\tilde{h}(t_1)=h(t_1,0,\ldots,0)$, and $\res_0$ denotes the
usual residue at $0$ of a meromorphic one-form in one variable.

The reason for bringing in the local trace is the following
commutative diagram \cite{L}:
$$
\begin{array}{rcl}
H^n_{\{p\}}(V,\Omega_V^n) & \stackrel{\cong}{\longrightarrow} &
H^n(V,\Omega_V^n)\\
 & & \\
tr_p\searrow\cong & & \cong\swarrow tr\\
 & & \\
 & \C &
\end{array}
$$
To tie in with the above description of $tr_p$, it remains to
explain the isomorphism
$$
H^n(V,\Omega_V^n)\longrightarrow H^n_{\{p\}}(V,\Omega_V^n)
\longrightarrow H^{n-1}(U-\{p\},\Omega_U^n)
$$ in \v{C}ech cohomology.
Take an affine (and hence acyclic) open covering
$\U=\{U_{\alpha}\}$ of $V$, so that
$H^n(V,\Omega_V^n)\cong\check{H}^n(\U,\Omega_V^n)$.
We may assume that the neighborhood $U$ of $p$ selected above is
entirely contained in some $U_{\alpha}$.
Then, throwing in $U_i=U-D_i$ and $U_0:=U$, we get an acyclic refinement
 $\bar{\U}$ of $\U$. Hence
$\check{H}^n(\bar{\U},\Omega_V^n)$ is still isomorphic to
$H^n(V,\Omega_V^n)$.

Now, an $(n-1)$-cochain over $U-\{p\}$ with coefficients in
$\Omega_U^n$ is a section of $\Omega_U^n$ over $U_1\cap\ldots\cap
U_n$. Regarding it as a section of $\Omega_V^n$ over $U_0\cap
U_1\cap\ldots\cap U_n$ defines an $n$-cochain with coefficients in
$\Omega_V^n$. This is the cochain map underlying the isomorphism
$$
H^{n-1}(U-\{p\},\Omega_U^n)\stackrel{\cong}{\longrightarrow}
H^n(V,\Omega_V^n)\ .
$$
To define its inverse, just read this backwards: take a cocycle in
$\check{C}^n(\bar{\U},\Omega_V^n)$, isolate its component in
$\Gamma(U_0\cap U_1\cap\ldots\cap U_n,\Omega_V^n)$, and view it as
an element in $\Gamma(U_1\cap\ldots\cap U_n,\Omega_U^n)$.

\section{The Laurent data}
Let $p$ be a regular point on $D$, and $\X$ a germ of a
smooth curve through $p$, normal to $D$, cut out by $t_2=\ldots=t_n=0$.
$z:=t_1|_{\X}$ defines a
coordinate on $\X$, identifying $\hat{\OO}_{\X,p}$ with $\C[[z]]$. We will
use the notation $\H=\C((z))$, $\H_+= \C[[z]]$, and  $\H'=\H/\C$,
and we will write $u$
for the Laurent expansion in terms of $z$ of the global meromorphic
objects on $V$ {\em restricted to} $\X$. Thus we have two
compositions, both denoted $u$:
$$
\Gamma(V,\OO_V(*D))\surj
\OO_{\X,p}(*p)\hookrightarrow\H\surj \H'
$$
and
$$
\Gamma(V,\Omega_V(*D))\surj \Omega^1_{\X,p}(*p)\hookrightarrow\H dz\ .
$$

Now we introduce two important subspaces of
$\H'$ carrying information about $V$:
$K_0:=u(\Gamma(V,\OO_V(*D)))$
and
$$
\Omega:=\{f\in\H'\,|\,df\in u(\Gamma(V,\tilde{\Omega}^1_V(*D)))\}\
{}.
$$
It is obvious that $K_0\subset\Omega$ and that $u$ (followed by
$d^{-1}:\H dz\stackrel{\sim}{\rightarrow}\H'$) induces a surjection
\begin{equation}
\label{basic:surj}
H^1(V,\C)\cong\frac{\Gamma(V,\tilde{\Omega}^1_V(*D))}
{d\Gamma(V,\OO_V(*D))}\surj\frac{\Omega}{K_0}\ .
\end{equation}
We also have surjective morphisms
$$
H^0(V,\Omega_V^1)=\Gamma(V,\tilde{\Omega}^1_V)\longrightarrow
\Omega\cap\H'_+
$$
and
$$
H^1(V,\OO_V)\cong\frac{\Gamma(V,\tilde{\Omega}^1_V(*D))}
{\Gamma(V,\Omega^1_V)+d\Gamma(V,\OO_V(*D))}\longrightarrow
\frac{\Omega}{(\Omega\cap\H'_+)+K_0}\ .
$$

The space $\H$ carries a symplectic form
$$
<f,g>=\res_0\,(fdg)\ ,
$$
which is non-degenerate on $\H'$. This form induces a symplectic
pairing, also denoted $<\ ,\ >$:
\begin{equation}
\label{form:symp}
(\Omega\cap\H'_+)\times
\frac{\Omega}{\Omega\cap\H'_+}\longrightarrow\C\ .
\end{equation}
A comparison of this with the polarization form $Q$ (\ref{Q})
constitutes our main result.

\section{The generalized reciprocity law}
\begin{Thm}
\label{main}
Let $\omega\in\Gamma(V,\Omega_V^1)$ and
$\eta\in\Gamma(V,\tilde{\Omega}_V^1(*D))$ represent classes in
$H^0(V,\Omega_V^1)$ and
$H^1(V,\OO_V)\cong \Gamma(V,\tilde{\Omega}_V^1(*D))/
\Gamma(V,\Omega_V^1)+d\Gamma(V,\OO_V^1(*D))$,
respectively, and suppose $u(\omega)=df$ and $u(\eta)=dg$ for some
$f\in\Omega\cap\H'_+$ and $g\in\Omega$. Then
$$
Q([\omega],[\eta])=(-1)^{n-1}2\pi i\,<f,g>\ .
$$
\end{Thm}

\pf
We will work in the covering $\bar{\U}$ as above.
The class $[\omega]$ is represented by the
cochain$\{\omega_{\alpha}=\omega|_{U_{\alpha}}\}$,
whereas $[\eta]=[\{g_{\alpha\beta}\in\Gamma(U_{\alpha}\cap
U_{\beta},\OO_V)\}]$ has been described as the end of Section
\ref{H1}. We single out
$$
g_{01}=\mu_1|_{U_0\cap U_1}-\mu_0|_{U_0\cap U_1}\ :
$$
since $U_0=U$ is
contractible, we may assume $d\mu_0=\eta|_{U}$, and since $\eta$ is
holomorphic over $U_1$, we may take $\mu_1=0$. Thus
$g_{01}=-\mu_0|_{U_0\cap U_1}$ is a holomorphic function on
$U_0\cap U_1$, with $\eta|_{U_0\cap U_1}=-dg_{01}$.

And $c_{[D]}\in H^1(V,\Omega_V^1)$ is represented by
\begin{equation}
\label{c}
\left\{-\frac{1}{2\pi i}
\frac{dh_{\alpha\beta}}{h_{\alpha\beta}}+d\ell\right\}\ ,
\end{equation}
where $h_{\alpha\beta}=0$ is a local equation of $D$ in
$U_{\alpha}\cap U_{\beta}$ and $\ell\in\Gamma(U_{\alpha}\cap
U_{\beta},\OO^*_V)$ extends across $D$ as an invertible holomorphic
function, if $D$ meets the closure of $U_{\alpha}\cap U_{\beta}$.

As explained in Section \ref{polarization},
$$
Q([\omega],[\eta])=(2\pi i)^n tr\,([\omega]\smile[\eta]\smile
c_{[D]}^{n-1})\ ,
$$and the trace can be computed locally, by taking the residue of
the component at
$U_0\cap U_1\cap\ldots\cap U_n$
of a \v{C}ech cochain representing
$[\omega]\smile[\eta]\smile c_{[D]}^{n-1}$ in $H^n(V,\Omega_V^n)$.
Thus it suffices to take the wedge product of the restrictions to
$U_0\cap U_1\cap\ldots\cap U_n$ of $\omega_0$, $g_{01}$, and of the
\v{C}ech components of $c_{[D]}$ over $U_0\cap U_1$, $U_1\cap U_2$,
\ldots,$U_{n-1}\cap U_n$.
However, in the case of the cochain (\ref{c}) we have no
information on the singularities these components have along the
$D_j$'s, except the one over $U_1\cap U_2$, which is
$-\frac{1}{2\pi i}\frac{dt_1}{t_1}+d\ell$.
This is unsuitable for computing the residue. Thus we need other
cochains representing $c_{[D]}$.
Now, replacing a divisor by a linearly equivalent one does not
affect the Chern class.
And just as $v_1t_1^N\in\OO_{V,p}$ is a germ of a global
meromorphic function on $V$ for an appropriate $N$ and some
$v_1\in\Gamma(U,\OO^*_U)$,
we shall assume temporarily that the
same is true of the other coordinate functions $t_2,\ldots,t_n$.
Then $v_{j1}(t_j/t_1)^N$ ($j=2,\ldots,n$) also come from meromorphic
functions on $V$.
Consequently, for each
$j=2,\ldots,n$ there is a cochain representing $c_{[ND]}$ whose
component on $U_{j-1}\cap U_j$
is $-\frac{N}{2\pi i}\frac{dt_j}{t_j}+d\ell_j$, with $\ell_j$
extending as an invertible holomorphic function on all of $U$.
Dividing by $N$ gives new cochains representing $c_{[D]}$.
Using these cochains,
\begin{eqnarray}
\nonumber
\lefteqn{Q([\omega],[\eta])=}\\
\nonumber & =  &
(2\pi i)^n \Res\ \omega_0 g_{01}\wedge
\left(-\frac{1}{2\pi
i}\frac{dt_2}{t_2}+d\ell_2\right)\wedge\ldots\wedge
\left(-\frac{1}{2\pi i}\frac{dt_n}{t_n}+d\ell_n\right)\\
 & = & (-1)^{n-1}2\pi i\, \Res\ \omega_0 g_{01}
\wedge\frac{dt_2}{t_2}\wedge\ldots\wedge\frac{dt_n}{t_n}
\label{Q:Res}
\end{eqnarray}
The last reduction is possible because $d\ell_j$'s do not
contribute to the residue, having no poles or zeroes along any
component of $D^+$.

At this point we recall that the residue is independent of the
parameter system (see, e.g. \cite{L}). Thus we may return to the
one in which $t_1$ is a local equation of $D$, while
$t_2,\ldots,t_n$ are arbitrary.

By (\ref{Res:res}) the multidimensional residue in (\ref{Q:Res})
reduces to the one-dimen\-sion\-al
$$
(-1)^{n-1}2\pi i\,
\res_0\,(g_{01}\omega_0|_{\X})=
(-1)^{n-1}2\pi i\,<-g,f>\ ,
$$
and we end up with
$$
Q([\omega],[\eta])=(-1)^{n-1}2\pi i\,<f,g>\ .
$$
\qed

\begin{Cor}
\label{zero}
 $<x,y>=0$ whenever $x\in\Omega\cap \H'_+$ and $y\in K_0$.
Consequently, the pairing (\ref{form:symp}) induces
$$
(\Omega\cap\H'_+)\times
\frac{\Omega}{\Omega\cap\H'_++K_0}\longrightarrow\C\ ,
$$
also denoted $<\ ,\ >$, and the surjection
$$
H^0(V,\Omega_V^1)\times H^1(V,\OO_V)\longrightarrow
(\Omega\cap\H'_+)\times \frac{\Omega}{\Omega\cap\H'_++K_0}
$$
transforms the pairing $Q$ on the source into $(-1)^{n-1}2\pi i\,<\
,\ >$ on
the target.
\end{Cor}
\qed

\begin{Cor}
\label{main:2}
The surjection
$$
H^1(V,\C)\cong\frac{\Gamma(V,\tilde{\Omega}_V^1(*D))}
{d\Gamma(V,\OO_V(*D))}\longrightarrow\frac{\Omega}{K_0}
$$
is a symplectic isomorphism.
\end{Cor}

\pf
The map in question induces a graded surjection
\begin{eqnarray*}
\lefteqn{H^0(V,\Omega_V^1)\oplus H^1(V,\OO_V)\cong}\\
& \cong &\Gamma(V,\Omega_V^1)\oplus
\frac{\Gamma(V,\tilde{\Omega}_V^1(*D))}
{\Gamma(V,\Omega_V^1)+d\Gamma(V,\OO_V(*D))}\\
  & \longrightarrow  &
(\Omega\cap\H'_+)\oplus\frac{\Omega}{\Omega\cap\H'_++K_0}
\end{eqnarray*}
which transforms the polarization form $Q$ on the left into the
residue pairing $(-1)^{n-1}2\pi i\,<\ ,\ >$ on the right. However,
$Q$ is non-degenerate, i. e. for any non-zero
$x\in H^0(V,\Omega_V^1)$ there exists $y\in H^1(V,\OO_V)$ with
$Q(x,y)\neq 0$, and similarly for any non-zero $y\in H^1(V,\OO_V)$.
Therefore, the image of any such $x$ or $y$ in $\Omega\cap\H'_+$
(resp. in $\Omega/\Omega\cap\H'_++K_0$) cannot be $0$. This shows
that the graded map is an isomorphism --- in fact, a symplectic
one, --- which implies the same for the original map. \qed

\section{Extended Hodge structure}
In this last section we will show that the above constructions can
be completed to a full {\em extended Hodge structure of weight one}
(= {\em an extended abelian variety}).

We already have $K_0\subset\H'$. Now we use the isomorphism
$H^1(V,\C)\cong\Omega/K_0$ to define $U$ as the image of
$H^{0,1}\subset H^1(V,\C)$ in $\Omega/K_0$. Put $Z=$ the preimage
of $U$ in $\Omega$ under the projection
$\Omega\rightarrow\Omega/K_0$.
Evidently, $U$ is a complement of $\Omega\cap\H'_+$ in
$\Omega/K_0$, and $Z$ is a complement of $\H'_+$ in $\H'$. And
since $U$ projects isomorphically onto $\Omega/\Omega\cap\H'_++K_0$
under the projection
$$
\frac{\Omega}{K_0}\longrightarrow\frac{\Omega}{\Omega\cap\H'_++K_
0}\ ,
$$
$U$ is perfectly paired with $\Omega\cap\H'_+$ under the pairing
induced by $<\ ,\ >$ on $\Omega/K_0$. From this we deduce that $Z$
is a maximal isotropic subspace of $\H'$. Furthermore,
$\Omega\subseteq K_0^{\perp}$, i. e. $Z\subset K_0^{\perp}$. Thus
$Z$ is sandwiched between $K_0$ and $K_0^{\perp}$:
$$
K_0\subset Z \subset K_0^{\perp}\ .
$$

\begin{Lemma}
The subspaces $K_0$ and $\Omega$ of $\H'$ are annihilators of each
other with respect to the symplectic form $<\ ,\ >$.
\end{Lemma}
\pf
By (\ref{zero}) $K_0\subseteq \Omega^{\perp}$. And by
(\ref{main:2}) and by non-degeneracy of $Q$,  ${<\ ,\ >}$ induces
a
nondegenerate pairing on $\Omega/K_0$. However,
$<\ ,\ >$ also induces a nondegenerate pairing on
$\Omega/\Omega^{\perp}$, which is a surjective image of
$\Omega/K_0$. But any symplectic surjection of  a vector space with
a non-degenerate symplectic form must be an isomorphism, as the
proof of (\ref{main:2}) shows. Hence $\Omega^{\perp}=K_0$.

Now, $\Omega\subseteq K_0^{\perp}$, hence
$K_0^{\perp\perp}\subseteq\Omega^{\perp}=K_0$. But $K_0\subseteq
K_0^{\perp\perp}$; so $K_0=K_0^{\perp\perp}$. Then $<\ ,\ >$ is
nondegenerate on $K_0^{\perp}/K_0^{\perp\perp}=K_0^{\perp}/K_0$. In
particular, its maximal isotropic subspaces must be of dimension
$\frac{1}{2}\dim\,(K_0^{\perp}/K_0)$.
However, $Z/K_0$ is already a maximal isotropic subspace in
$K_0^{\perp}/K_0$, and it is only of dimension
$g=\frac{1}{2}\dim\,(\Omega/K_0)$. Therefore,
$K_0^{\perp}=\Omega$. \qed

Combining this lemma with Corollary \ref{main:2} yields
\begin{Thm}
\label{main:two}
There is a symplectic isomorphism
$$
(H^1(V,\C),\ Q)\longrightarrow
(K_0^{\perp}/K_0,\ (-1)^{n-1}2\pi i<\ ,\ >)\ .
$$
\end{Thm}

We now  introduce the remaining components of the extended abelian
variety. Define $\Lambda\subset\Omega/K_0=K_0^{\perp}/K_0$ as the
image of the lattice $H^1(V,\Z)$ under the isomorphism
$H^1(V,\C)\stackrel{\sim}{\rightarrow}K_0^{\perp}/K_0$, and let $K$
be the preimage of $\Lambda$ under the projection
$K_0^{\perp}\rightarrow K_0^{\perp}/K_0$.

With this notation we may summarize our results as follows.
\begin{Thm}
\label{EHS}
 The triple
$(Z,K_0,\Lambda)$ associated to the pointed polarized variety
$(V,D,p)$ is an extended abelian variety.
\end{Thm}
Indeed, the definition of an extended abelian variety in \cite{AD}
calls
for $Z$ to be a maximal
isotropic subspace of $\H'$ with $\H'_+\oplus Z=\H'$ , $K_0$ must
be a subspace of $Z$ and $\Lambda$ a lattice in $K_0^{\perp}/K_0$
such that
$$
(\Lambda,K_0^{\perp}/K_0,2\pi i\,<\ ,\ >)
$$
constitutes
a polarized Hodge structure of weight one, with the Hodge
decomposition induced by the direct sum decomposition of $\H'$:
the $(1,0)$-component of $K_0^{\perp}/K_0$ is
$K_0^{\perp}\cap\H'_+$, and the $(0,1)$-component is
${(K_0^{\perp}\cap Z)/K_0}$. These conditions have been established
already.

\rk
In \cite{SW} and \cite{AD} the polarization form $Q$ on the first
cohomology
of a Riemann surface corresponds to $\frac{1}{2\pi i}<\ ,\ >$
instead of our
$2\pi i\,<\ ,\ >$. The discrepancy is due to a different convention
adopted in
these papers: they identify $\Lambda$ with $H^1(X,2\pi i\Z)$ rather
than
$H^1(X,\Z)$, as we do.

\rk
\label{choices}
The construction of an extended abelian variety associated to $V$
obviously depends on the choice of the ample divisor $D$ and the
point
$p$. It also depends on the coordinate system $(t_1,\ldots,t_n)$ at
$p$,
without which
we would not be able to define the Laurent expansions and hence the
map $u$. However, all such choices are equally good for our
purposes.

We end with one last observation. In the curve case we had
$$
K=\{f\in\H'\mid e^{2\pi if}\in u(\Gamma(X-\{p\},\OO_X^*))\}=
u(\Gamma(X,\OO_X(*p)/\Z))\ .
$$
It turns out, $K$ admits a similar identification in the
multidimensional situation.

\begin{Lemma}
$K=u(\Gamma(V,\OO_V(*D)/\Z))$.
\end{Lemma}
\pf
Let us write $\tilde{K}$ for $u(\Gamma(V,\OO_V(*D)/\Z))$.
Our point of departure in identifying $K$ with $\tilde{K}$ is the
exact sequence
$$
0\longrightarrow\Z_V\longrightarrow\OO_V(*D)\longrightarrow\OO_V(
*D)/\Z\longrightarrow 0\ .
$$
The corresponding cohomology sequence reads, in part,
$$
 H^0(V,\OO_V(*D))\longrightarrow H^0(V,\OO_V(*D)/\Z)\longrightarrow
H^1(V,\Z)\longrightarrow H^1(V,\OO_V(*D))
$$
The last term is isomorphic to $H^1(V-D,\OO_V)=0$. Therefore,
applying $u$ yields
$$
0\longrightarrow K_0\longrightarrow \tilde{K}\longrightarrow
\tilde{K}/K_0\longrightarrow 0\ .
$$
Thus $u$ identifies $H^1(V,\Z)$ with $\tilde{K}/K_0$. It is easy to
see that $\tilde{K}\subset\Omega$, and we conclude that
$\tilde{K}=K$. \qed

\rk
Connecting $\Gamma(V,\OO_V(*D)/\Z)$ with $\Gamma(V-D,\OO_V^*)$ by
means of the  exponential sequence on $V-D$, we also get
$$
K=\{f\in\H'\mid e^{2\pi i f}\in u(\Gamma(V-D,\OO_V^*))\}\ .
$$

\

E-mail address: {\em yk@yu1.yu.edu}

\end{document}